\documentclass[12pt]{amsart}
\usepackage{amscd}
\usepackage{amsmath}

\usepackage{amssymb}

\numberwithin{equation}{section}

\newtheorem{thm}[equation]{Theorem}
\newtheorem{prop}[equation]{Proposition}

\newtheorem{lemma}[equation]{Lemma}
\newtheorem{cor}[equation]{Corollary}

\theoremstyle{definition}

\newcommand{\Gal}{\mathop{\mathrm{Gal}}}

\newcommand{\Z}{\mathbb{Z}}
\newcommand{\F}{\mathbb{F}}

\newcommand{\Spec}{\operatorname{Spec}}

\newcommand{\Aut}{\operatorname{Aut}}

\newcommand{\PS}{\mathbb{P}}

\newcommand{\Oplus}{\operatornamewithlimits{\textstyle\bigoplus}}

\renewcommand{\phi}{\varphi}

\DeclareMathAlphabet{\cat}{OT1}{cmss}{m}{sl}

\title
{Artin shapes}


\author
{Nikita Karpenko}

\address
{Mathematical \& Statistical Sciences \\
University of Alberta \\
Edmonton
\\
CANADA}

\email
{karpenko@ualberta.ca}
\urladdr{www.ualberta.ca/~karpenko}

\author
{Guangzhao Zhu}

\address
{Mathematical \& Statistical Sciences \\
University of Alberta \\
Edmonton
\\
CANADA}

\email
{guangzha@ualberta.ca}

\date
{18 Nov 2024}

\begin{document}

\begin{abstract}
We introduce and study on examples a notion of the Artin shape for a motive related to a projective homogenous
variety. We apply it to the problem of finding the complete motivic decomposition of the variety.
Our examples cover
unitary involution varieties as well as some varieties given by a quadratic Weil transfer.
Some of the decompositions obtained dispel prior expectations on how motivic decompositions
of projective homogeneous varieties can look like.
\end{abstract}

\maketitle

\tableofcontents

\addtocounter{section}{-1}

\section{The scope}
\label{The scope}

We consider Grothendieck's Chow motives (defined as in \cite[\S64]{EKM}) over a field
with coefficients in $\F:=\Z/p\Z$, where $p$ is  a prime number.
Following a tradition, we also write $\F$ for the motive of the spectrum of the base field and call a {\em Tate motive}
its (Tate) shift $\F\{i\}$ with any $i\in\Z$.
Another important type of objects are direct summand in the spectra of \'etale algebras, called
{\em Artin motives}.

Let $M$ be a direct summand in the motive of a
projective homogeneous variety under a reductive group $G$ over a field $F$
and let $B$ be the variety of Borel subgroups in $G$.
By \cite[Theorem 7.5]{MR2110630},
over the function field $F(B)$, $M$ becomes a direct sum of (Tate) shifts of indecomposable
Artin motives which we call the {\em Artin shape} (or simply the {\em shape}) of $M$.
As shown in \cite{MR2264459} (see also \cite[Corollary 2.6]{upper}),
the decomposition we use to define the Artin shape is unique.
Besides,
the motive $M$ itself decomposes in a finite direct sum of indecomposables and
such a decomposition, called {\em complete}, is also unique.

The shape of the entire motive $M(X)$ of a projective $G$-homogeneous variety $X$,
which we also call the {\em motivic shape} of $X$,
contains the Artin motive given by the Tate motive $\F:=\Spec(F)$ with no shift (i.e., with the shift $0$),
whereas all remaining Artin motives in the shape have positive shifts.
It follows that the complete motivic decomposition of $X$ contains precisely one summand $U(X)$,
called the {\em upper motive} of $X$,
whose shape satisfies the same property.

By \cite[Corollary 2.15]{upper},
the upper motives of projective homogeneous varieties satisfy a very efficient isomorphism criterion:
$U(X)\simeq U(X')$ if and only if each of the varieties $X_{F(X')}$ and $X'_{F(X)}$ possesses a
closed point of prime to $p$ degree.

A description of summands that can appear in the complete decomposition of $M(X)$,
already previously available in a different situation (cf.\! \cite[Theorem 1.1]{outer}), has been recently
provided in \cite{aum32a} for {\em $p'$-inner} $G$, i.e., for $G$ that acquires inner type over a finite Galois
field extension $E/F$ of prime to $p$ degree and whose {\em higher Tits $p$-indexes},
defined in \cite{titspindexes}, are stable
under the action of the Galois group of $E/F$ on the Dynkin diagram of $G$:

\begin{thm}[{\cite[Theorem 7.1 and Remark 6.9]{aum32a}}]
\label{Aupper}
Every summand in the complete motivic decomposition of a projective homogeneous variety $X$ under a $p'$-inner group $G$ is a shift of the tensor product $U(Y)\!\otimes\! A$, where $Y$ is a projective $G$-homogeneous variety with $X(F(Y))\ne\emptyset$ and $A$ is an Artin motive isomorphic to an indecomposable summand in the motive
of the spectrum of the $F$-algebra $L$ given by an intermediate field $L$ in $E/F$.
\end{thm}

In this note we work out several new applications for that result.

Assume that $p$ is odd and fix a quadratic Galois field extension $L/F$.
The motive of the spectrum of the $F$-algebra $L$ splits off the Tate motive $\F$;
we write $A$ for the complementary summand.
This $A$ is an indecomposable Artin motive non-isomorphic to $\F$ over $F$ but becoming isomorphic to
$\F$ over $L$.

For some $n\geq1$,
we let $D$ be a degree $p^n$ central division $L$-algebra, write $X$ for its Severi-Brauer
$L$-variety, and consider the $F$-variety $R(X)$ given by the Weil transfer of $X$.
The variety $R(X)$ is projective homogeneous under the reductive group $G=R(\Aut(D))$ given by the Weil
transfer of the automorphism group of $D$.
This $G$ acquires inner type over $L$.
Since $p$ is odd, the degree of the field extension $L/F$ is prime to $p$.
However, $G$ is not necessarily $p'$-inner:
it is $p'$-inner if and only if the algebra
$D$ is {\em balanced} in the sense of \cite[\S4]{zinovy-k} meaning that the algebra $\check{D}$ conjugate
to $D$ is Brauer-equivalent to a tensor power of $D$.
Note that if this is the $r$th tensor power, then $r^2\equiv1\mod{p^n}$ implying that $r\equiv\pm1\mod{p^n}$.
It follows that $\check{D}$ is isomorphic either to $D$ or to its opposite algebra.
(And the two conditions are mutually exclusive.)
Considering the norm algebra of $D$, defined in \cite[\S3B]{MR1632779},
one translates the first condition as $D\simeq D'_L$ for a degree $p^n$ central division $F$-algebra $D'$.
The second condition means that the norm algebra is Brauer-trivial and translates
by \cite[Theorem 3.1(2)]{MR1632779} as
``$D$ admits a unitary $F$-involution''.

We are going to investigate the shape of the {\em upper motive} $U(R(X))$
as well as the complete motivic decompositions of the variety $R(X)$ and of another related variety
in these two cases.
Our main results are Theorems \ref{main1}, \ref{main2}, and \ref{main3}.


\section{The first case}
\label{The first case}

Assume that we are in the first case, i.e., $D\simeq D'_L$,  and let $X'$
be the Severi-Brauer $F$-variety of $D'$.
Then $U(R(X))\simeq U(X')$ by the isomorphism criterion mentioned in \S\ref{The scope}.
By \cite[Corollary 2.22]{upper} (see also \cite[Theorem 2.2.1]{GrCh}),
the motive $M(X')$ is indecomposable.
Therefore $U(X')=M(X')$.
The variety $X'$ is projective homogeneous under the group $G':=\Aut(D')$.
Over the function of field of the corresponding variety of Borel subgroups $X'$ becomes isomorphic
to a ($p^n-1$)-dimensional projective space $\PS$ whose motivic shape is
$$
\F\oplus\F\{1\}\oplus\dots\oplus\F\{p^n-1\}.
$$
So, this is also the shape of the upper motive $U(R(X))$ we started with.
It can be represented ``graphically'' as
$$
\begin{array}{cccc}
\{0\} & \{1\} & \dots & \{p^n-1\}\\
\F & \F & \dots & \F
\end{array}
\hspace{15ex}\text{ or simply }\hspace{10ex}
\begin{array}{cccc}
 &  &  &\\
\F & \F & \dots & \F.
\end{array}
$$

The graphical shape of $M(R(X))$ is
\begin{equation}
\label{shapeR}
\begin{array}{ccccccccc}
\F &  \F  & \dots & \F  &     & &\\
    &  A   &  A  & \dots &  A & &\\
    &       &  \F &  \F  & \dots & \F  &     & &\\
    &       &       & A   &  A  & \dots &  A & &\\
    &       &       &      &   \ddots   &     &      &\ddots &\\
   & &       &       &       &  \F & \F  & \dots & \F
\end{array}
\end{equation}
(with $p^n$ variously shifted lines, each of length $p^n$).
Indeed, this shape is given by the complete motivic decomposition of the Weil transfer
$R(\PS)$.
The decomposition is determined by the action of the Galois group
$\Gal(L/F)$ on the Chow group of the product $R(\PS)_L=\PS\times\PS$, performing
via the factors exchange.

By Theorem \ref{Aupper}, the complete decomposition of $M(R(X))$ consists of shifts of
$U(R(X))$ and of $U(R(X))\!\otimes\! A$.
If follows
that the complete decomposition of $M(R(X))$ (over $F$) is
\begin{multline}
\label{M}
M(R(X))\simeq U(R(X))\;\oplus\;\big(U(R(X))\!\otimes\! A\big)\{1\}\;\oplus\;\\
U(R(X))\{2\}\;\oplus\;\big(U(R(X))\!\otimes\! A\big)\{3\}\;\oplus\dots\oplus
U(R(X))\{p^n-1\}\;=\\
\left(\Oplus_{i=0}^{(p^n-1)/2} U(R(X))\{2i\}\right)\;\;\oplus\;\;\left(\Oplus_{i=1}^{(p^n-1)/2}\big(U(R(X))\!\otimes\! A\big)\{2i-1\}\right).
\end{multline}

We put together the main result of this section dispelling some prior expectation on how motivic decompositions
of projective homogeneous varieties can look like:

\begin{thm}
\label{main1}
Let $p$ be an odd prime, $n$ a positive integer, $X'$ the Severi-Brauer variety of a degree $p^n$ central
division algebra over a field $F$.
For a quadratic Galois field extension $L/F$, formula
(\ref{M}) provides the complete motivic decomposition for the Weil transfer
$R(X)$ of the $L$-variety $X:=X'_L$.
The decomposition holds in the category of Chow motives with coefficients in $\Z/p\Z$.
\end{thm}

\section{The second case}
\label{The second case}

Now let us assume that $D$ admits a unitary $F$-involution.
First of all, let us note that the shape of $M(R(X))$ is still given by (\ref{M}):
the proof remains unchanged.
What we want to figure out is the shape of $U(R(X))$ and the complete decomposition of $M(R(X))$.
Concerning the latter, by Theorem \ref{Aupper} once again we still know  a priori that
$M(R(X))$ is a direct sum of shifts of $U(R(X))$ and of $U(R(X))\!\otimes\!A$.

By \cite[Lemma 6.4]{aum32a},
the motive $U(R(X))$ remains indecomposable over $L$ so that $U(R(X))_L\simeq U(R(X)_L)$.
Since $U(R(X)_L)\simeq U(X)=M(X)$ has the shape $\F\;\F\;\dots\;\F$,
we conclude that
the shape of $U(R(X))$ looks as $\F**\dots*$, where each of $*$ is either $\F$ or $A$.
We will be able to precisely determine the shape
in the next section based on the following preliminary result:

\begin{lemma}
\label{lemma}
The shape of $U(R(X))$ is either $\F\;\F\;\dots\;\F$ (as in \S\ref{The first case})
or $\F\;A\;\,\F\;A\;\dots\;\F\;A\;\,\F$
($(p^n+1)/2$ exemplars of $\F$ alternated with $(p^n-1)/2$ exemplars of $A$).
The complete decomposition of  $M(R(X))$ in the first case is (\ref{M});
in the second case it is given by
\begin{multline}
\label{second}
M(R(X))\simeq \\
U(R(X))\;\oplus\;U(R(X))\{1\}\;\oplus\;\dots\;\oplus\;U(R(X))\{p^n-1\}=
\Oplus_{i=0}^{p^n-1}U(R(X))\{i\}.
\end{multline}
\end{lemma}

\begin{proof}
Assuming that the shape of $U:=U(R(X))$ is $\F\;\F**\dots*$, i.e., that the component with shift $\{1\}$ is given by
the Tate motive,
we prove by induction on $i\geq3$ and looking at (\ref{shapeR})
that the component with shift $\{i\}$ is also $\F$.
At the same time, we prove the claim on the complete decomposition of $M(R(X))$.
As the initial step, we see that the complete decomposition starts with
$U\;\oplus\; (U\!\otimes\!A)\{1\}$ implying that $U$ has the shape $\F\;\F\;\F**\dots*$.
Knowing this, we get  that the next term of the complete decomposition is $U\{2\}$ and therefore the shape of
$U$ has $\F$ in position $\{3\}$.
And so on.

Similarly, assuming now the shape of $U$ starts with $\F\;A$,
we prove by induction
that the shape has to be $\F\;A\;\,\F\;A\;\dots\;\F$ and the complete motivic decomposition of $R(X)$
is given by (\ref{second}).
\end{proof}

At this point, we do not know yet if decomposition (\ref{second}) actually occurs, but, as already
mentioned, we will prove that it does in the next section.
Let us illustrate this result by noticing that
decompositions (\ref{M}) and (\ref{second}) correspond to the following two different graphical
representation of the motivic shape of $R(X)$:
$$
\begin{array}{ccccccccc}
\F &  \F  & \dots & \F  &     & &\\
    &  A   &  A  & \dots &  A & &\\
    &       &  \F &  \F  & \dots & \F  &     & &\\
    &       &       & A   &  A  & \dots &  A & &\\
    &       &       &      &   \ddots   &     &      &\ddots &\\
   & &       &       &       &  \F & \F  & \dots & \F
\end{array}
\hspace{1em}=\hspace{1em}
\begin{array}{ccccccccc}
\F &  A  & \F & A & \dots & \F  &     & &\\
    &  \F   &  A  & \F& A &\dots &  \F & &\\
       &      &   \ddots   &     &      & &  &\ddots& \\
   & &      &  \F & A  & \F & A & \dots & \F
\end{array}
$$

Let us also note that the ratio
\begin{equation}
\label{ratio}
\frac{\text{number of $\F$}}{\text{number of $A$}}
\end{equation}
in the second possible shape of $U(R(X))$ equals $(p^n+1)/(p^n-1)$ and is the same as in the shape of $M(R(X))$.
Note that the similar ratio for $U(Y)\!\otimes A\!$ equals $(p^n-1)/(p^n+1)$ which is lower.
This gives another proof of (\ref{second}) in the second case of Lemma \ref{lemma}:
we already know that $M(R(X))$ is a direct sum of shifts of $U(R(X))$ and of
$U(R(X))\!\otimes\!A$, but due to the ratio considerations, the tensor product with $A$ cannot appear.

\section{The unitary involution variety}

In the situation of \S\ref{The second case}, let us fix a unitary $F$-involution $\tau$ on the $L$-algebra $D$.
The reductive $F$-group $G:=\Aut_L(D,\tau)$
(called the projective unitary group of $(D,\tau)$)
is $p'$-inner and acquires inner type over $L$.
Let us determine the complete motivic decomposition of the
corresponding {\em involution variety} $Y$ -- the variety of
rank $1$ (over $L$) right isotropic ideals in $D$.
We start with some information on its motivic shape:

\begin{prop}
\label{propo}
Given any $i=1,2,\dots,(p^n-1)/2$,
the motivic shape of the variety $Y_i$ of rank $i$ (over $L$) right isotropic ideals in the algebra $D$
contains
$(b_i+a_i)/2$ shifts of the Tate motive $\F$ and of $(b_i-a_i)/2$ shifts of the Artin motive $A$, where
$$
a_i=2^i\cdot\binom{(p^n-1)/2}{i}\;\;\text{ and }\;\;
b_i=\binom{p^n}{i}\cdot\binom{p^n-i}{i}.
$$
\end{prop}

\begin{proof}
Here is the meaning of the integers $a_i$ and $b_i$ appearing in the statement of Proposition \ref{propo}.
The integer $b_i$ is the {\em rank} of $M(Y_i)$, i.e., the total number of motives in its shape.
To get the formula for $b_i$,
note that over a separable closure of the base field, $Y_i$ becomes the variety of flags in a $p^n$-dimensional vector space, where the flags are given by an $i$-plane contained in a ($p^n-i$)-plane.

To explain the meaning of $a_i$, we first note that over $F(B)$ the motive $M(Y_i)$ is also a direct sum of shifts
of the Tate motive $\F$
and of the motive of spectrum of the $F(B)$-algebra $L(B)$; the integer $a_i$ is the number of shifts of $\F$.
It follows from \cite[Lemma 7.1]{gug}
for unitary grassmannians and the similar formula \cite[Formula 2.6]{sgog} for orthogonal grassmannians
that $a_i$ coincides with the rank of the motive of an orthogonal grassmannian of totally isotropic $i$-planes in a $p^n$-dimensional non-degenerate quadratic form, indicated, e.g., in
\cite[Proof of Theorem 2.2]{Kresch} (see also \cite[Proof of Theorem 1.2]{Kresch}).
\end{proof}

\begin{cor}
\label{cor}
The number of $\F$ in the shape is $(p^n+1)(p^n-1)$ whereas the number of
$A$ is $(p^n-1)^2$.
In particular,  ratio (\ref{ratio}) for the motivic shape of $Y=Y_1$ is
equal to  $(p^n+1)/(p^n-1)$.
\qed
\end{cor}

Note that $U(Y)$ is isomorphic to the upper motive $U(R(X))$.
Therefore, the two possible shapes of $U(Y)$ are listed in Lemma \ref{lemma}.

\begin{thm}
\label{main3}
Let $p$ be an odd prime, $n$ a positive integer,
$L/F$ a quadratic Galois field extension,
$D$ a degree $p^n$ central division $L$-algebra endowed with a unitary $F$-involution $\tau$,
and $Y$ the corresponding involution variety (over $F$).
Then the shape of $U(Y)$ is as in the second case of Lemma \ref{lemma} whereas
the complete motivic decomposition of $Y$ is
\begin{equation}
\label{third}
M(Y)\simeq
\Oplus_{i=0}^{p^n-2}U(Y)\{i\}.
\end{equation}
\end{thm}

In contrast to (\ref{second}), $i$ in (\ref{third}) ranges up to $p^n-2$ only.

\begin{proof}[Proof of Theorem \ref{main3}]
By Theorem \ref{Aupper}, every summand in the complete motivic decomposition of $Y$ is a shift of
$U(Y)$ or of $U(Y)\!\otimes\!A$.
Therefore,
if $U(Y)$ has the first shape of Lemma \ref{lemma},
the number of $\F$ in the shape of $M(Y)$ is divisible by $p^n$.
On the other hand, it is not by Corollary \ref{cor}.
This contradiction proves the statement of Theorem \ref{main3} on the shape of $U(Y)$.

Now, by the same argument as in \S\ref{The second case}, we get that the complete decomposition
of $M(Y)$ consists of shifts of $U(Y)$ only.
Over $L$, the variety $Y$ becomes a rank $p^n-2$ projective bundle over the Severi-Brauer variety $X$
and so the motive of $Y$ decomposes over $L$ as in (\ref{third}).
It follows that the same decomposition holds already over $F$.
\end{proof}

For the illustrating purpose, let us draw the shape of the motive of $Y$
the way reflecting its complete decomposition:
$$
\begin{array}{ccccccccc}
\F &  A  & \F & A & \dots & \F  &     & &\\
    &  \F   &  A  & \F& A &\dots &  \F & &\\
       &      &   \ddots   &     &      & &  &\ddots& \\
   & &      &  \F & A  & \F & A & \dots & \F
\end{array}
$$
($p^n-1$ variously shifted lines, each of length $p^n$).

With Theorem \ref{main3}, we also proved our last main result:

\begin{thm}
\label{main2}
In the situation of Theorem \ref{main3}, let $X$ be the Severi-Brauer variety of $D$ and let $R(X)$
be its Weil transfer to $F$.
Then its upper motive $U(R(X))$ has the second shape of Lemma \ref{lemma} whereas its total motive
$M(R(X))$ has the complete decomposition (\ref{second}).
\qed
\end{thm}

Decompositions of Theorem \ref{main2} and \ref{main3} put together yield
$$
M(R(X))\simeq M(Y)\oplus U(Y)\{p^n-1\}.
$$


\begin{thebibliography}{10}

\bibitem{Kresch}
{\sc Buch, A.~S., Kresch, A., and Tamvakis, H.}
\newblock Quantum {P}ieri rules for isotropic {G}rassmannians.
\newblock {\em Invent. Math. 178}, 2 (2009), 345--405.

\bibitem{MR2110630}
{\sc Chernousov, V., Gille, S., and Merkurjev, A.}
\newblock Motivic decomposition of isotropic projective homogeneous varieties.
\newblock {\em Duke Math. J. 126}, 1 (2005), 137--159.

\bibitem{MR2264459}
{\sc Chernousov, V., and Merkurjev, A.}
\newblock Motivic decomposition of projective homogeneous varieties and the
  {K}rull-{S}chmidt theorem.
\newblock {\em Transform. Groups 11}, 3 (2006), 371--386.

\bibitem{titspindexes}
{\sc De~Clercq, C., and Garibaldi, S.}
\newblock Tits {{\(p\)}}-indexes of semisimple algebraic groups.
\newblock {\em J. Lond. Math. Soc., II. Ser. 95}, 2 (2017), 567--585.

\bibitem{aum32a}
{\sc De~Clercq, C., Karpenko, N., and Qu\'eguiner-Mathieu, A.}
\newblock A-upper motives of reductive groups.
\newblock Preprint (1 Nov 2024, 23 pages). Available on the 2nd author's
  webpage.

\bibitem{EKM}
{\sc Elman, R., Karpenko, N., and Merkurjev, A.}
\newblock {\em The algebraic and geometric theory of quadratic forms}, vol.~56
  of {\em American Mathematical Society Colloquium Publications}.
\newblock American Mathematical Society, Providence, RI, 2008.

\bibitem{GrCh}
{\sc Karpenko, N.~A.}
\newblock Grothendieck {C}how motives of {S}everi-{B}rauer varieties.
\newblock {\em St. Petersburg Math. J. 7}, 4 (1996), 649--661.

\bibitem{outer}
{\sc Karpenko, N.~A.}
\newblock Upper motives of outer algebraic groups.
\newblock In {\em Quadratic forms, linear algebraic groups, and cohomology},
  vol.~18 of {\em Dev. Math.} Springer, New York, 2010, pp.~249--258.

\bibitem{sgog}
{\sc Karpenko, N.~A.}
\newblock Sufficiently generic orthogonal {G}rassmannians.
\newblock {\em J. Algebra 372\/} (2012), 365--375.

\bibitem{gug}
{\sc Karpenko, N.~A.}
\newblock Unitary {G}rassmannians.
\newblock {\em J. Pure Appl. Algebra 216}, 12 (2012), 2586--2600.

\bibitem{upper}
{\sc Karpenko, N.~A.}
\newblock Upper motives of algebraic groups and incompressibility of
  {S}everi-{B}rauer varieties.
\newblock {\em J. Reine Angew. Math. 677\/} (2013), 179--198.

\bibitem{zinovy-k}
{\sc Karpenko, N.~A., and Reichstein, Z.}
\newblock A numerical invariant for linear representations of finite groups.
\newblock {\em Comment. Math. Helv. 90}, 3 (2015), 667--701.
\newblock With an appendix by Julia Pevtsova and Reichstein.

\bibitem{MR1632779}
{\sc Knus, M.-A., Merkurjev, A., Rost, M., and Tignol, J.-P.}
\newblock {\em The book of involutions}, vol.~44 of {\em American Mathematical
  Society Colloquium Publications}.
\newblock American Mathematical Society, Providence, RI, 1998.
\newblock With a preface in French by J.\ Tits.

\end{thebibliography}

\end{document}